\documentclass[11pt]{article}
\usepackage{amsmath}
\usepackage{amstext, latexsym, amssymb}

\newcounter{question}[section]
\newcommand \R{\mathbb{R}}
\DeclareMathOperator{\Vol}{Vol} 

\newcommand{\cone}{\mbox{$\times \hspace*{-0.228cm} \times$}}

\begin{document}
\baselineskip 7.2 truemm

\title{\textbf {\Large\center DENSITY OF A MINIMAL SUBMANIFOLD AND
\\
TOTAL CURVATURE OF ITS BOUNDARY}}

\normalsize
\author{\bf Jaigyoung Choe\thanks{\noindent The first-named author
was supported in part by KRF-2007-313-C00057.}\, and Robert
Gulliver}
\date{}
\maketitle \baselineskip=1.2\normalbaselineskip \noindent


\begin{abstract}

Given a piecewise smooth submanifold $\Gamma^{n-1} \subset \R^m$
and $p \in \R^m$, we define the {\em vision angle} $\Pi_p(\Gamma)$
to be the $(n-1)$-dimensional volume of the radial projection of
$\Gamma$ to the unit sphere centered at $p$. If $p$ is a point on
a stationary $n$-rectifiable set $\Sigma \subset \R^m$ with
boundary $\Gamma$, then we show the density of $\Sigma$ at $p$ is
$\leq$ the density at its vertex $p$ of the cone over $\Gamma$. It
follows that if $\Pi_p(\Gamma)$ is less than 
twice
the volume of $S^{n-1}$, for all $p \in \Gamma$, then $\Sigma$ is
an embedded submanifold.  As a consequence, we prove that given
two $n$-planes $R^n_1, R^n_2$ in $\R^m$ and two compact convex
hypersurfaces $\Gamma_i$ of $R^n_i, i=1,2$, a nonflat minimal
submanifold spanned by $\Gamma:=\Gamma_1\cup\Gamma_2$ is
embedded.\\

\end{abstract}


\section{Introduction}

Fenchel \cite{F1} showed that the total curvature of a closed
space curve $\gamma\subset\R^m$ is at least $2\pi$, and it
equals $2\pi$ if and only if $\gamma$ is a plane convex curve.
F\'{a}ry \cite{Fa} and Milnor \cite{M} independently proved that a
simple knotted regular curve has total curvature larger than
$4\pi$. These two results indicate that a Jordan curve which is
curved at most {\it double} the minimum
is isotopically simple. But in fact minimal
surfaces spanning such Jordan curves must be simple as well.
Indeed, Nitsche \cite{N} showed that an analytic Jordan curve in
$\R^3$ with total curvature at most $4\pi$ bounds exactly
one minimal disk. Moreover, Ekholm, White and Wienholtz \cite{EWW}
proved that a minimal surface spanning such a Jordan curve in
$\R^m$ is embedded.

Given an $n$-dimensional submanifold $M$ of $\R^m$, there
are two well-studied ways of defining the total curvature of $M$:
the higher-dimensional Gauss-Bonnet integral $\int_M\Omega$ as
defined in \cite{AW} and \cite{C1}; and the total
absolute curvature of $M$, $\int_M K^*dV_M$ as defined by
Chern and Lashof in \cite{CL} (see section 2 below).
Chern and Lashof proved that
$\int_{M}K^*\, dV_M\geq 2, $ with equality if and only if $M$ is a
convex hypersurface in an $(n+1)$-dimensional plane.  Eells and
Kuiper have shown
that if $\int_{M}K^*dV_M<3$ then $M$ is homeomorphic to
$\mathbb{S}^n$ and that if $\int_MK^*dV_M<4$ then $M$ is 
homeomorphic to $\mathbb S^n$, $\R P^n$, $\mathbb CP^{n/2},$
$\mathbb HP^{n/4}$ or to $CayP^2$ (for $n=16$). 
\cite{EK}.  

In the light of Ekholm-White-Wienholtz's theorem, it is quite
natural to conjecture that {\it an $n$-dimensional minimal
submanifold $\Sigma\subset\R^m$ spanning a compact
connected submanifold $\Gamma^{n-1}$ with total absolute curvature
$<4$ is embedded.} In this paper we prove a theorem in the spirit
of this conjecture: given two $n$-planes $R^n_1,R^n_2$ in $\R^m$
and two compact convex hypersurfaces $\Gamma_i^{n-1}$ of $R^n_i,
i=1,2$, a nonflat minimal submanifold spanned by
$\Gamma:=\Gamma_1\cup\Gamma_2$ is 
embedded.

In \cite{Fa} F\'{a}ry showed that the total curvature of a space
curve $\gamma$ in $\R^m$ is equal to the average over all
2-planes $R^2\subset\R^m$ of the total curvature of the
orthogonal
projection of $\gamma$ onto the $R^2$. We shall use an extension
of F\'{a}ry's theorem, due to Langevin and Shifrin \cite{LS},
which shows that given an $(n-1)$-dimensional
submanifold $\Gamma$ of $\R^m$, the total absolute curvature
of $\Gamma$ equals the average over all $n$-planes
$R^n\subset\R^m$ of the total absolute curvature of the orthogonal
projection of $\Gamma$ onto the $n$-plane $R^n$.

\section{Total absolute curvature}\label{sec2}
Consider a submanifold $M^n$ of Euclidean space $\R^m$.
As discussed above, in high dimension and codimension we discuss
two types of total curvature:  one
intrinsic (Allend\"orfer-Weil-Chern-Gauss-Bonnet),
and one extrinsic (Chern-Lashof). In this
section we shall review Chern-Lashof's total absolute curvature.
This total curvature may be understood in terms of Gauss-Kronecker
curvature of hypersurfaces.

Let $M^n$ be an oriented hypersurface immersed in
$\R^{n+1}$. A unit normal vector $\nu$ to $M$ at $p\in M$
defines the Gauss map $G_1:M\rightarrow\mathbb{S}^{n}$. The
determinant of the differential $G_{1}*$, or of the second
fundamental form of $M$, is called the
{\it Gauss-Kronecker curvature} of $M$, which we shall denote
$GK_M$. It follows that for $M$ compact,
$$\int_MGK_M\,dV_M=c_n\,{\rm deg}(G_1),~~c_n:=
{\Vol}(\mathbb{S}^n).$$
Furthermore, if $n$ is even, H. Hopf \cite{H} showed
\begin{equation}\label{Hopf}
\int_MGK_M\,dV_M=\frac{1}{2}\,c_n\,\chi(M).
\end{equation}

Now let $M$ be an $n$-dimensional submanifold of $\R^m$. The
volume form of the unit normal bundle $N_1M$ of $M$ is $dV_M\wedge
d\sigma_{m-n-1}$ where the restriction of $d\sigma_{m-n-1}$ to a
fiber of $N_1M$ at $p$ is the volume form of the sphere of unit
normal vectors at $p\in M$. Define the Gauss map
$G_1:N_1M\rightarrow \mathbb{S}^{m-1}$ by $G_1(p,\nu)=\nu$ and let
$d\sigma_{m-1}$ be the volume form of $\mathbb S^{m-1}$. Then the
{\it Lipschitz-Killing curvature} $G(p,\nu)$ of $M$ at $(p,\nu)$
is defined to be the scalar $G(p,\nu)$ such that
$$G_1^*(d\sigma_{m-1})=G(p,\nu)\, dV_M\wedge d\sigma_{m-n-1}.$$
Then $G(p,\nu)$ is exactly
the volume expansion ratio of $G_1$, that is,
$$G(p,\nu)=\lim_{D\rightarrow \{p\}}\frac{{\rm Vol}(G_1(D))}{{\rm
Vol}(D)},$$ where ${\rm Vol}(G_1(D))$ denotes the signed volume of
$G_1(D)$. In fact, $G(p,\nu)$ has the following geometric
interpretation \cite{CL}: $G(p,\nu)$ is equal to the
Gauss-Kronecker curvature at $p$ of the orthogonal projection of
$M$ onto the $(n+1)$-dimensional plane $L(\nu)$ spanned by $T_pM$
and $\nu$.

Let $\pi$ be the canonical projection of $N_1M$ into $M$. The
integrals
$$K(p):=\frac{1}{c_{m-1}}\int_{\pi^{-1}(p)}
G(p,\nu)\, d\sigma_{m-n-1}~~~{\rm and}~~~
K^*(p):=\frac{1}{c_{m-1}}\int_{\pi^{-1}(p)} |G(p,\nu)|\,
d\sigma_{m-n-1}$$ are called the {\it total curvature} and the
{\it total absolute curvature} of $M$ at $p$, respectively. The
integrals
$$\tau(M):=\int_{M} K\, dV_M,~~~{\rm and}~~~
\tau^*(M):=\int_{M} K^*\, dV_M$$ are called the {\it total
curvature} and the {\it total absolute curvature of $M$},
respectively. Lipschitz and Killing have shown that $K(p)$ is an
intrinsic quantity of $M$ at $p$ for $n$ even (see \cite{SS} for a
more general result). However, $K(p)=0$ for $n$ odd. Both
$\tau(M)$ and $  \tau^*(M)$ remain unchanged even if the ambient
space $\R^m$ is embedded into $\R^k,~ k>m$.

For $M^n\subset\R^m$, Fenchel \cite{F2} generalized Hopf's
theorem \eqref{Hopf}:
\begin{equation}\label{Fenchel}
\int_M K \, dV_M=\chi(M).
\end{equation}
In contrast, Chern and Lashof \cite{CL} proved that
\begin{equation}\label{C} \int_{M}K^*\, dV_M\geq 2, \end{equation}
with equality if and only if $M$ is a convex hypersurface in an
$(n+1)$-dimensional plane, and that if $\int_{M}K^*dV_M<3$ then
$M$ is homeomorphic to $\mathbb{S}^n$. Moreover,
Morse theory tells us that
$$\int_M K^*\,dV_M\geq \sum_i\beta_i,$$
where $\beta_i$ is the $i$-th Betti number of $M$ ([W], Theorem
28).

\section{Vision angle versus average density}\label{sec3}

A minimal submanifold $\Sigma^n$ in $\R^m$ has the remarkable
property that the density  of $\Sigma$ at $p\in\Sigma$ is bounded
above by that of the cone $C=p\cone\partial\Sigma$ at its vertex
$p$. (We assume that $\Sigma$ with its boundary is compact.)
Recall that the {\it density} of $\Sigma$ is defined as
$$ \Theta_\Sigma(p) = \lim_{r \to 0} \frac{\Vol\Big(\Sigma \cap
B_r^m(p)\Big)}{\Vol\Big(B_r^n(p)\Big)}. $$
Further, the density of a cone $C$ has the interesting property
that it
equals the average of the densities of the orthogonal projections
of $C$ onto $n$-planes in $\R^m$. These
properties will be verified in this section.\\

In what follows, we shall write ${\overline\nabla}$ for the
Euclidean connection on $\R^m$, and $\nabla = \nabla_M$ for the
induced connection on a submanifold $M$.\\

\noindent{\bf Lemma 1.} {\it Let $\Sigma$ be an $n$-dimensional
minimal submanifold of $\R^m$, $p$ a point of
$\R^m$, and $C$ an $n$-dimensional
piecewise smooth 
cone with vertex $p$.
Define the Euclidean distance function
$r(x)={\rm dist}(p,x),x\in\R^m$. Let
$Y_1=r\overline{\nabla}r$ and
$Y_2=r^{1-n}\overline{\nabla}r$, and define
${\rm div}_\Sigma
Y_i={\rm tr}_\Sigma\overline{\nabla}Y_i=\sum_j\langle
\overline{\nabla}_{e_j}Y_i,e_j\rangle$, $\{e_1,\dots,e_n\}$ being
an
orthonormal frame of $\Sigma$.   Then

\noindent{\rm (a)} On $\Sigma$, ${\rm div}_\Sigma Y_1=n$ and ${\rm
div}_\Sigma Y_2\geq 0$;

\noindent{\rm (b)} On $C$, ${\rm div}_C Y_1=n$ and ${\rm div}_C
Y_2=0$.}\\

We require that $C$ be piecewise smooth, that is, a topological
manifold which has a triangulation into simplices that are $C^2$ up
to their boundaries.\\

\noindent{\it Proof.} Given an $n$-dimensional submanifold
$M\subset \R^m$, it is well known that
$$\triangle_M \,x:=(\triangle_M\, x_1,\dots,\triangle_M\,
x_m)=\vec{H},$$
where $\vec{H}$ is the mean curvature vector of $M$,
the trace of its second fundamental form.
Hence the orthogonal coordinate functions $x_1,\dots,x_m$ of
$\R^m$ are
harmonic on a minimal submanifold $\Sigma^n$ of $\R^m$. If
we take $p$ as the origin, then since $\vec{H} = 0$ on $\Sigma,$
$$\,{\rm div}_\Sigma(Y_1)={\rm div}_\Sigma(r \overline{\nabla}r)=
\frac{1}{2}\triangle_\Sigma\, r^2+\langle
r\overline{\nabla}r,\vec{H}\rangle=\frac{1}{2}\sum\triangle_\Sigma
\,x_i^2=\sum x_i\triangle_\Sigma\, x_i+\sum|\nabla x_i|^2=n.$$
On the cone $C$, since $\vec{H}$ is perpendicular to
$r\overline{\nabla}r=x\in C$, we have
$$ \,{\rm div}_C(Y_1)={\rm div}_C (r\overline{\nabla}r)=
\frac12\triangle_C\, r^2+\langle r\overline{\nabla}r,\vec{H}\rangle
=\frac12\sum\triangle_C \,x_i^2=
\langle x, \vec{H}\rangle+\sum|\nabla x_i|^2=n.$$
On the other hand, for 
$M=\Sigma$ or $C$,
\begin{eqnarray*}
{\rm div}_M\,Y_2&=& {\rm div}_M (r^{-n}Y_1)=
-n r^{-n-1}\langle \nabla r, Y_1\rangle + r^{-n} {\rm div}_M (Y_1)
= nr^{-n}\Big( -|\nabla r|^2+1\Big).
\end{eqnarray*}

Note that $|\nabla r|\leq 1$ on $M=\Sigma$ and $|\nabla r|\equiv
1$ on $M=C$. This completes the proof.~~~~~~~~$\square$\\

\noindent{\bf Theorem 1.} {\it Let $\Sigma$ be a 
stationary $n$-rectifiable set  with boundary $\Gamma$
in $\R^m$, an
open dense subset of $\Sigma$ being a smooth minimal submanifold.
Let $C$ be the cone $p \cone\Gamma,\, p\in\R^m$. Then
$$\Theta_\Sigma(p)\leq \Theta_C(p),$$
with equality if and only if $\Sigma=C$ and $C$ is star-shaped
with respect to $p$.}\\

\noindent{\it Proof.} Compute the first variation of volume with
respect to the (Lipschitz continuous) variation vector field
$$ Y: = r^{1-n}\overline{\nabla} r ~~~~{\rm for}~~ r \geq
\varepsilon $$
 and
$$ Y: = \varepsilon^{-n} r \overline{\nabla} r~~~~
{\rm for}~~ r \leq \varepsilon.  $$
Then the first variation of $\Sigma$ with respect to the flow with
velocity field $Y$ [Si, p. 80] is
$$ \int_\Sigma {\rm div_\Sigma} Y \, dV_\Sigma, $$
which must equal
$$ \int_{\Gamma} \langle Y, \nu_\Sigma\rangle \, dV_\Gamma, $$
where $\nu_\Sigma$ is the outward unit normal vector to $\Gamma$
tangent to $\Sigma$.

  Computing the divergence on smooth subsets of the stationary set
$\Sigma$, we find by Lemma 1 (a)
\begin{equation}\label{inequality}
{\rm div_\Sigma} Y \geq 0~~~~{\rm for}~~r \geq
\varepsilon,\end{equation}
with equality at points where $\overline{\nabla}r$ lies in the
tangent space, and
$$
{\rm div_\Sigma}Y = n \varepsilon^{-n}~~~~{\rm for}~~r \leq
\varepsilon.
$$
It follows that for each small $\varepsilon$,
\begin{equation}\label{d1}
\frac{{\rm Vol}(\Sigma\cap B_\varepsilon(p))}{|B^n_1|
\varepsilon^n} \leq \frac{1}{n|B^n_1|} \int_\Gamma
r^{1-n}\langle\overline{\nabla}r,\nu_\Sigma\rangle\,dV_\Gamma,
~|B_1^n|:={\rm Vol}(B^n_1(0)).
\end{equation}
Now apply Stokes' theorem to the integral of ${\rm div}_CY$
on $C$:
$$\int_C{\rm div}_CY\,dV_C=\int_{\partial C}\langle
Y,\nu_C\rangle=\int_\Gamma\langle Y,\nu_C\rangle,$$
where $\nu_C$ is the outward unit conormal to $\Gamma$ on $C$.
Therefore, by Lemma 1 (b)
\begin{equation}\label{d2}
\frac{{\rm Vol}(C\cap B_\varepsilon(p))}{|B_1^n| \varepsilon^n}
= \frac{1}{n|B_1^n|} \int_\Gamma
r^{1-n}\langle\overline{\nabla}r,\nu_C\rangle\,dV_\Gamma.
\end{equation}
Note here that
$$0\leq\langle\overline{\nabla}r,\nu_C\rangle$$
and
\begin{equation}\label{inequality2}
\langle \overline{\nabla}r,\nu_\Sigma\rangle\leq
\langle\overline{\nabla}r,\nu_C\rangle.\end{equation}
Thus, letting $\varepsilon\rightarrow 0$ in inequality \eqref{d1}
and
equation \eqref{d2}, we get the desired density estimate. If
equality
holds, then we must have equality in inequalities
\eqref{inequality} and
\eqref{inequality2}, which implies $\Sigma=C$ and $\partial
r/\partial\nu\geq0$.~~~~~~~~$\square$
\\

\noindent {\bf Definition 1.} Let $\pi_p$ be the radial projection
of $\R^m\setminus\{p\}$ onto $\partial B_1(p)$, the unit
sphere centered at $p\in\R^m$. Define the {\it vision
angle at $p$}
of an $(n-1)$-rectifiable set $\Gamma\subset\R^m$ by
$$\Pi_p(\Gamma)={\rm Vol}(\pi_p(\Gamma)),$$
and the {\it vision angle} of $\Gamma$ by
$$\Pi(\Gamma)={\rm sup}_{\,p\in\R^m}\Pi_p(\Gamma)).$$
Here the volume ${\rm Vol}(\pi_p(\Gamma))$ counts multiplicity.\\

Clearly we have for any $p\in\R^m$ and $C:=p\cone\Gamma$
$$c_{n-1}\Theta_C(p)=
\Pi_p(\Gamma^{n-1})\leq\Pi(\Gamma), ~~~c_{n-1}:= {\rm
Vol}(\mathbb{S}^{n-1}),$$
and hence we get the following corollaries to Theorem 1.\\

\noindent{\bf Corollary 1.} {\it If $\Gamma\subset\R^m$ is
an $(n-1)$-dimensional compact manifold, then any
stationary rectifiable set
$\Sigma$ spanning $\Gamma$ satisfies}
$$c_{n-1}\Theta_\Sigma(p)\leq\Pi_p(\Gamma)$$
{\it for all $p \in \Sigma.$}\\

\noindent{\bf Corollary 2.} {\it If $\Gamma\subset\R^m$ is
an $(n-1)$-dimensional compact manifold with $\Pi(\Gamma) < 2
c_{n-1}$, then any immersed minimal submanifold
$\Sigma$ spanning $\Gamma$ is embedded.}\\

\noindent{\it Proof}. An immersed submanifold $\Sigma$ with
density $\Theta_\Sigma(q)<2$ at each point $q \in \R^m$
has no self-intersection.~~~~~~~~$\square$\\

\noindent{\bf Remark.} It may appear inappropriate
to view $\Pi(\Gamma)$ as a total curvature. But it has its own
merit, as the following example demonstrates. Define an immersed
closed $C^1$ curve $\gamma\subset \R^2$ (the unit square
plus four small loops at the corners) by        
$$\gamma=
\partial([-1,1]^2)\cup\{(x,y):|x|>1, |y|>1, [(|x|-1)^2 +
(|y|-1)^2]^{3/2}=
\varepsilon (|x|-1)(|y|-1)\}$$
%
and define a Jordan curve $\Gamma\subset\R^n$ to be
an embedded  $C^2$ curve $C^1$-close    
to $\gamma$. Then for small $\varepsilon,$
$$\int_\Gamma|\vec{k}|ds >6\pi,~~~~  
{\rm however,}~~~~\Pi(\Gamma)\approx3\pi.$$
Hence by Corollary 2 any immersed 
minimal surface $\Sigma$ spanning $\Gamma$ is embedded
since $2c_1 = 4\pi$, although the  
Ekholm-White-Wienholtz theorem \cite{EWW} cannot give the same
conclusion.\\

Let $G_n(\R^m)$ denote the {\it Grassmann manifold} of
$n$-planes through the origin in $\R^m$, equipped with the
unique ${\mathbb O}(m)$-invariant probability measure, 
and let ${\rm
Ave}_{P\in G_n(\R^m)}$ be the average over all $P\in
G_n(\R^m)$. Denote by $\psi_P$ the orthogonal projection
of $\R^m$ onto $P\in G_n(\R^m)$.\\

\noindent{\bf Lemma 2.} {\it Let $\mathbb{S}^{n-1}$ be the unit
sphere in $\R^n\subset\R^m$ centered at the origin
$O$ of $\R^m$ and let $D$ be a domain in
$\mathbb{S}^{n-1}$. Then
$${\rm Ave}_{P\in G_n(\R^m)}\{\Theta_{\psi_P(O\cone
D)}(O)\}=\Theta_{O\cone D}(O).$$}

\noindent{\it Proof.} Assume that $a(D)>0$ is a positive real
number
such that
\begin{equation}\label{avedensity}
{\rm Ave}_{P\in G_n(\R^m)}\{\Theta_{\psi_P(O\cone
D)}(O)\}=a(D)\cdot\Theta_{O\cone D}(O).
\end{equation}
Letting $D$ shrink to a point $x\in\mathbb{S}^{n-1}$, one can
define a function $a:\mathbb{S}^{n-1}\rightarrow\R$ given
by
$$a(x):=\lim_{D\rightarrow\{x\}}\frac{{\rm
Ave}_{P\in G_n(\R^m)}\{\Theta_{\psi_P(O\cone
D)}(O)\}}{\Theta_{O\cone D}(O)}.$$
Then, by means of a partition of unity by functions of small
support, one can see that
$$a(D)=\frac{\int_{D}a(x)dV_{\mathbb{S}^{n-1}}}{{\rm Vol}(D)}.$$
Note here that $\mathbb{O}(n)$ is transitive on $\mathbb{S}^{n-1}$
and that the elements of $\mathbb{O}(n)$ preserve the volume form
$dV_{\mathbb{S}^{n-1}}$ on $\mathbb{S}^{n-1}$. Therefore one
concludes that for all $x \in \mathbb{S}^{n-1}$,
$$a(x)\equiv c~~~{\rm for~a~positive~constant~ }c$$
 and hence for any domain $D\subset\mathbb{S}^{n-1}$,
 $$a(D)\equiv c.$$
Therefore it follows from equation \eqref{avedensity} that
\begin{equation}\label{avec}
{\rm Ave}_{P\in G_n(\R^m)}\{\Theta_{\psi_P(O\cone
D)}(O)\}=c\cdot\Theta_{O\cone D}(O)
\end{equation}
for any domain $D\subset\mathbb{S}^{n-1}$. However, for almost all
$P\in G_n(\R^m)$,
$$\Theta_{\psi_P(O\cone\mathbb{S}^{n-1})}(O)=
\Theta_{O\cone\mathbb{S }^{n-1}}(O)=1.$$
Thus $c=1$ in equation \eqref{avec}, which  completes the
proof.~~~~~~~~$\square$\\

\noindent{\bf Theorem 2.} {\it Let $\Gamma^{n-1}\subset \R^m$
be a compact submanifold. Then}
$${\it \Pi_q(\Gamma^{n-1})=
{\rm Ave}_{P\in
G_n(\R^m)}\,\{\Pi_{\psi_P(q)}(\psi_P(\Gamma))\}.}$$

\noindent{\it Proof.} The cone $q\cone\Gamma$ 
can be thought of as a union of infinitesimal cones
$q\cone\Delta\Gamma_i$ and then one can apply Lemma 2 to each
$q\cone\Delta\Gamma_i$. Hence
\begin{eqnarray*}
\Pi_q(\Gamma)&=&c_{n-1}\,\Theta_{q\cone\Gamma}(q)\\
&=&c_{n-1}\,{\rm Ave}_{P\in
G_n(\R^m)}\,\{\Theta_{\psi_P(q\cone\Gamma)}(\psi_P(q))\}\\
&=&{\rm Ave}_{P\in
G_n(\R^m)}\,\{\Pi_{\psi_P(q)}(\psi_P(\Gamma))\}.~~~~~~~~\square
\end{eqnarray*}


We shall also require the following generalization of F\'{a}ry's
theorem to any dimension $n$ and to any codimension $m-n$, which
was proved by Langevin and Shifrin ([LS], Proposition 2.15): \\

\noindent{\bf Theorem LS.} {\it Let $\Gamma^{n-1}$ be a smooth
submanifold of $\R^m, m\geq n$. Then}
$$\frac{c_{n-1}}{2}\int_\Gamma K^*dV_\Gamma=
{\rm Ave}_{P\in
G_n(\R^m)}\int_{\psi_P(\Gamma)}|GK_{\psi_P(\Gamma)}|
\,dV_{\psi_P(\Gamma)}.$$\\

\section{Embeddedness of minimal submanifolds}\label{sec4}

It is tempting to propose a higher-dimensional extension of
Ekholm-White-Wienholtz's theorem as follows:\\

\noindent {\bf Conjecture.}
If $q \in \Sigma,$  a minimal submanifold of $\R^m$ spanning an
$(n-1)$-dimensional compact manifold $\Gamma$, then
$$\Theta_\Sigma(q)\leq\frac12\int_\Gamma K^*dV_\Gamma.$$

If this were known, one could prove the following as well:\\

{\it If an $(n-1)$-dimensional compact connected manifold $\Gamma$
satisfies $\int_\Gamma K^*dV_\Gamma<4$, then any immersed 
minimal submanifold $\Sigma^n$ spanning $\Gamma$ is embedded.}\\

Conjecture seems to be hard to prove as yet. 
 
However, if we let
$\Gamma_i$ be a compact convex hypersurface of an affine $n$-plane
$R^n_i\subset\R^m$, $i=1,2$, and define
$\Gamma=\Gamma_1\cup\Gamma_2$, then we may prove Conjecture for
this case. Our proof uses the vision angle of $\Gamma$ from a
point of $\Sigma$, and averages over projections onto all
$n$-dimensional subspaces $P$ of $\R^m$.
Namely, for $i=1,2$,
\begin{equation}\label{flat}
\Pi_{\psi_P(q)}(\psi_P(\Gamma_i))\leq
c_{n-1}=\int_{\psi_P(\Gamma_i)}|GK_{\psi_P(\Gamma_i)}|dV_{\psi_P(\Gamma_i)},
\end{equation}
since $\psi_P(\Gamma_i)$ is a convex hypersurface in
$\psi_P(R^n_i)$.  Here, equality holds for all $P$ if and only if
$q$ is in $R^n_i$ and inside $\Gamma_i$. Thus we have the
following:\\

\noindent {\bf Theorem 3.} {\it Given two $n$-planes $R^n_1,
R^n_2$ in $\R^m$, let $\Gamma_i$ be a compact convex
hypersurface in $R^n_i$, $i=1,2$. If
$\Gamma=\Gamma_1\cup\Gamma_2,$
then any $n$-dimensional minimal submanifold $\Sigma$ spanning
$\Gamma$ is either a union of two flat domains of $R^n_i$ or is
nonflat and has no self intersection.}\\

\noindent{\it Proof.} We may compute that
$\int_\Gamma K^*dV_\Gamma=
\sum_{i=1,2} \int_{\Gamma_i} K^* \, dV_{\Gamma_i} = 4$.
Thus by inequality \eqref{flat} and Corollary 1 we have
$\Theta_\Sigma\leq 2.$ If $\Theta_\Sigma=2$, inequality
\eqref{flat} and Corollary 1 imply $\Sigma$ is flat.
If $\Theta_\Sigma<2$ everywhere on $\Sigma$, then $\Sigma$
is nonflat and has no self intersection. ~~~~~~~~$\square$\\

\noindent {\bf Remark.} It should be mentioned that R. Schoen
\cite{Sc} proved a theorem
which implies a special case of Theorem 3:\\
{\it If $\Gamma=\Gamma_1\cup \Gamma_2$ where $\Gamma_1, \Gamma_2$
are $(n-1)$-spheres in parallel $n$-planes with the line
$\ell$ joining their centers being orthogonal to these
hyperplanes, then any immersed minimal submanifold $\Sigma^n$
spanning $\Gamma$ is a hypersurface of revolution with axis
$\ell$. In particular, $\Sigma$ is a catenoid or a pair of plane
disks.}

\vspace{1cm}
\small
\begin{tabbing}
aaaaaaaaaaaaaaaaaaaaaaaaaaaaaasssssssssssssssss \=
bbbbbbbbbbbbbbbbbbbbbbbbbbbbbbb\kill

Jaigyoung Choe                       \> Robert Gulliver \\
Korea Institute for Advanced Study   \> School of Mathematics \\
Hoegiro 87, Dongdaemun-gu            \> 127 Vincent Hall \\
Seoul, 130-722                    \> University of Minnesota\\
Korea                                \> Minneapolis MN 55414, USA
\\
email: {\tt choe@kias.re.kr}     \> {\tt gulliver@math.umn.edu}\\
fax: +82-2-958-3786                  \> +1-612-626-2017\\
{\tt http://newton.kias.re.kr/\~{ }choe}\>
                     {\tt http://www.ima.umn.edu/\~{ }gulliver} \\
\end{tabbing}

\end{document}